\documentclass[final,1p,times]{elsarticle}

\usepackage{amssymb}
 \usepackage{amsthm}
\usepackage{amscd}
\usepackage{amsmath}
\usepackage{amsfonts}
\usepackage{amssymb}
\usepackage{graphicx}
\usepackage{color}
\newtheorem{theorem}{Theorem}

\newtheorem{lemma}[theorem]{Lemma}
\newtheorem{claim}{Claim}
\newtheorem{remark}{Remark}

\usepackage{mathrsfs}
\usepackage{titletoc}

\newcommand{\ra}{\rightarrow}

\newcommand{\f}{\frac}

\renewcommand{\f}{\frac}

\newcommand{\be}{\begin{equation}}
\newcommand{\ee}{\end{equation}}
\newcommand{\bea}{\begin{eqnarray}}
\newcommand{\eea}{\end{eqnarray}}
\newcommand{\bna}{\begin{eqnarray*}}
\newcommand{\ena}{\end{eqnarray*}}

\renewcommand{\le}{\left}
\newcommand{\ri}{\right}

\journal{***}

\begin{document}

\begin{frontmatter}

\title{Multiple solutions of Kazdan-Warner equation on graphs in the negative case}

\author{Shuang Liu}
 \ead{shuangliu@ruc.edu.cn}
 \author{Yunyan Yang\footnote{Corresponding author}}
 \ead{yunyanyang@ruc.edu.cn}

 \address{Department of Mathematics,
Renmin University of China, Beijing 100872, P. R. China}

\begin{abstract}
Let $G=(V,E)$ be a finite connected graph, and let $\kappa: V\rightarrow \mathbb{R}$ be a function such that $\int_V\kappa d\mu<0$. We consider the following Kazdan-Warner equation on $G$:
\[\Delta u+\kappa-K_\lambda e^{2u}=0,\]
where $K_\lambda=K+\lambda$ and $K: V\rightarrow \mathbb{R}$ is a non-constant function satisfying $\max_{x\in V}K(x)=0$ and $\lambda\in \mathbb{R}$. By a variational method, we prove that there exists a $\lambda^*>0$ such that
when $\lambda\in(-\infty,\lambda^*]$ the above equation has solutions, and has no solution when $\lambda\geq \lambda^\ast$. In particular, it has only one solution if $\lambda\leq 0$; at least two distinct solutions if
$0<\lambda<\lambda^*$; at least one solution if $\lambda=\lambda^\ast$. This result complements earlier work of
Grigor'yan-Lin-Yang \cite{GLY16}, and is viewed as a discrete analog of that of Ding-Liu \cite{DL95} and Yang-Zhu \cite{YZ19} on manifolds.
\end{abstract}


\end{frontmatter}

\section{Introduction}
Variational method is always a powerful tool in partial differential equations and geometric analysis. Recently,
using this tool, Grigor'yan-Lin-Yang \cite{GLY16,GLY16',GLY17} obtained existence results for solutions to
various partial differential equations on graphs. In particular, Kazdan-warnar equation was proposed on graphs
in \cite{GLY16}. The Kazdan-Warner equation arises from the basic geometric problem on prescribing Gaussian curvature of Riemann surface, which systematically studied by Kazdan-Warner \cite{KW74,KW740}. On a closed Riemann surface
$(\Sigma,g)$ with the Gaussian curvature $\kappa$, let $\widetilde{g}=e^{2u}g$ be a smooth metric conformal to $g$
and $K$ be the Gaussian curvature
with respect to $\widetilde{g}$. Then $u$ satisfies the equation
\begin{equation}\label{eq:kw0}
\Delta_g u+\kappa-K e^{2u}=0,
\end{equation}
where $\Delta_g$ denotes the Laplace-Beltrami operator with respect to the metric $g$.
Let $v$ be a solution to $\Delta_gv=\overline{\kappa}-\kappa$ and $f=2(u-v)$, where $\overline{\kappa}$ is the averaged integral of
$\kappa$. Then the above equation is transformed to
\[\Delta_g f+2\overline{\kappa}-(2Ke^{2v})e^{f}=0.\]
Hence, one can free \eqref{eq:kw0} from the geometric situation, and just studies the equation
\begin{equation}\label{eq:kw1}
\Delta_g f+c-h e^{f}=0,
\end{equation}
where $c$ is a constant and $h$ is a function.
On graphs, it seems to be out of reach to resemble this topic in terms of Gaussian curvature. Therefore, in \cite{GLY16}, the authors focused on the equation similar to the form of \eqref{eq:kw1}, namely the Kazdan-Warner equation on graph, and obtained the following:
when $c=0$, it has a solution if and only if $h$ changes sign and the integral of $h$
is negative; when $c>0$, it  has a solution if and only if $h$ is positive somewhere; when
$c < 0$,  there is a threshold $c_h < 0$ such
that it has a solution if $c\in(c_h, 0)$, but it has no solution for any $c<c_h$. Later, Ge \cite{G17} found a solution in the critical case $c = c_h$. More recently Ge-Jiang \cite{G-J} studied the Kazdan-Warner equation on infinite graphs and Keller-Schwarz \cite{K-S} on
canonically compactifiable graphs; Camilli-Marchi \cite{C-M} extended the Kazdan-Warner equation on network; for other related works, we refer the readers to \cite{H-S-Z,L-W}.

Let us come back to a closed Riemann surface $(\Sigma,g)$, whose Euler characteristic is negative, or equivalently
$\int_\Sigma \kappa dv_g<0$. Replacing $K$ by $K+\lambda$ in (\ref{eq:kw0}) with
$K\leq 0$, $K\not\equiv 0$, and $\lambda\in\mathbb{R}$, Ding-Liu \cite{DL95} obtained the following conclusion by using a method of upper and lower solutions and a variational method: there exists a $\lambda^\ast>0$ such that if $\lambda\leq 0$, then (\ref{eq:kw0}) has a unique solution; if $0<\lambda<\lambda^\ast$, then (\ref{eq:kw0}) has at least two distinct solutions; if $\lambda=\lambda^\ast$, then  (\ref{eq:kw0}) has at least one solution; if $\lambda>\lambda^\ast$, then (\ref{eq:kw0}) has no solution. Recently, this result was partly reproved by
Borer-Galimberti-Struwe \cite{B-G-S} via a monotonicity technique due to Struwe \cite{Struwe1,Struwe2}, and was extended to the case of conical metrics by Yang-Zhu \cite{YZ19}.

 Our aim is to extend results of Ding-Liu \cite{DL95} to graphs. Let us recall some notations from graph theory.
Throughout this paper, $G=(V,E)$ is assumed to be a finite connected graph.  The edges on the graph are allowed to be weighted.  Weights are given by a function $\omega: V\times V\rightarrow [0,\infty)$, the edge $xy$ from $x$ to $y$ has weight $\omega_{xy}>0$.  We assume this weight function is symmetric, $\omega_{xy}=\omega_{yx}$.
Let $\mu:V \rightarrow \mathbb{R^{+}}$ be a positive measure on the vertices of the $G$. Denote by $V^{\mathbb{R}}$ the space of real functions on $V$. and by $\ell^{p}_\mu=\{f \in V^{\mathbb{R}}:\sum_{x\in V}\mu(x)|f(x)|^{p}<\infty\}$, for any $1\leq p< \infty$,
the space of $\ell^{p}$ integrable functions on $V$ with respect to the measure $\mu$.  For $p=\infty$, let $\ell^{\infty}=\{f \in V^{\mathbb{R}}:\sup_{x\in V}|f(x)|<\infty\}$ be the set of all bounded functions.
As usual, we define the $\ell^p_\mu$ norm of $f\in \ell^{p}_\mu,1\leq p\leq \infty$, by
\[\|f\|_p=\left(\sum_{x\in V}\mu(x)|f(x)|^p\right)^{{1}/{p}}, 1\leq p<\infty,\, \|f\|_{\infty}=\sup_{x\in V}|f(x)|.\]
We define the Laplacian $\Delta:V^{\mathbb{R}}\rightarrow V^{\mathbb{R}}$ on $G$ by
\be\label{Laplace}\Delta f(x)=\frac{1}{\mu(x)}\sum_{y\sim x} \omega_{xy}(f(x)-f(y)).\ee
Given the weight $\omega$ on $E$, there are two typical choices of Laplacian as follows:
\begin{itemize}
  \item $\mu(x)=\deg(x):=\sum_{y\sim x}\omega_{xy}$ for all $x\in V$, which is called the normalized graph Laplacian;
  \item $\mu(x)\equiv 1$ for all $x\in V$, which is the combinatorial graph Laplacian.
\end{itemize}
In this paper, we do not restrict $\mu(x)$ to the above two forms, but only require $\mu(x)>0$ for all $x\in V$.
Note that the Laplace operator defined in (\ref{Laplace}) is the negative usual Laplace operator.
The gradient form is defined by
\begin{eqnarray*}
2\Gamma(f,g)(x)& =&(f\cdot \Delta g+g\cdot\Delta f-\Delta(f\cdot g))(x)\\
               & =&\frac{1}{\mu(x)}\sum_{y\sim x}\omega_{xy}(f(x)-f(y))(g(x)-g(y)).
\end{eqnarray*}
For the sake of simplicity, we write $\Gamma(f,f)=\Gamma(f)$. Sometimes we use the notation $\nabla f \nabla g=\Gamma(f,g)$.  The length  of the gradient is denoted by
\[|\nabla f|(x)=\sqrt{\Gamma(f)}(x).\]
From now on, we write $\int_V ud\mu=\sum_{x\in V}\mu(x)u(x)$.
Define a Sobolev space with a norm on the graph $G$ by
\[W^{1,2}(V)=\left\{u\in V^{\mathbb{R}}:\int_V(|\nabla u|^2+u^2)d\mu<+\infty\right\},\]
and
\[\|u\|_{W^{1,2}(V)}=\left(\int_V(|\nabla u|^2+u^2)d\mu\right)^{1/2}\]
respectively. Since $G$ is a finite graph, we have that $W^{1,2}(V)$ is exactly $V^{\mathbb{R}}$, a finite dimensional linear space. This implies the following Sobolev embedding:
\begin{lemma}[\cite{GLY16}, Lemma 5]\label{lemma:compact}
If $G$ is a finite graph, then the Sobolev space $W^{1,2}(V)$ is precompact. Namely, if $\{u_j\}$ is bounded in $W^{1,2}(V)$, then there exists some $u\in W^{1,2}(V)$ such that up to a subsequence, $u_j\rightarrow u$ in $W^{1,2}(V)$.
\end{lemma}

The Kazdan-Warner equation we are interested in this paper reads as
\begin{equation}\label{eq:0}
\Delta u+\kappa-K_\lambda e^{2u}=0 \quad\mbox{on}\quad V,
\end{equation}
where $\kappa\in V^{\mathbb{R}}$ is a function, and $K_\lambda=K+\lambda$, $\lambda\in \mathbb{R}$,  $K\in V^{\mathbb{R}}$ is a function. Now we are ready to state our main results.
\begin{theorem}\label{thm:main}
Let $G=(V,E)$ be a finite graph, $\kappa$ and $K_\lambda$ be given as in \eqref{eq:0} such that
 $\int_V\kappa d\mu<0$, $K\leq \max_{ V}K=0$, and $K\not\equiv 0$.
 Then there exists a $\lambda^*\in (0,-\min_{V}K)$ satisfying
\begin{enumerate}
  \item if $\lambda\leq0$, then \eqref{eq:0} has a unique solution;
  \item if $0<\lambda<\lambda^*$, then \eqref{eq:0} has at least two distinct solutions;
  \item if $\lambda=\lambda^*$, then \eqref{eq:0} has at least one solution;
 \item if $\lambda>\lambda^*$, then \eqref{eq:0} has no solution;
\end{enumerate}
\end{theorem}
\begin{remark}
The assertion of $\lambda^*<-\min_{V}K$ comes from the conclusion of Step 2 in Subsection 3.3.
\end{remark}
\begin{remark}
Compared to the existence of solutions in the literature (see for example \cite{GLY16,G17}), the above results firstly reveal the multiple solution problem of Kazdan-Warner equation on graphs in the negative case.
\end{remark}
The proof of Theorem \ref{thm:main} is based on the method of variation. It can be viewed as a discrete analog
of the result of Ding-Liu \cite{DL95}.
The remaining part of this paper will be organized as follows: In Section \ref{sec-2}, we give several preliminary lemmas
for our use later; In Section \ref{sec-3}, we finish the proof of Theorem \ref{thm:main}.
\section{Preliminaries}\label{sec-2}
In this section, we provide discrete versions of the maximum principle, the Palais-Smale condition and the upper and lower solution principle. Note that $G=(V,E)$ is a finite connected graph.
\subsection{maximum principle}
To proceed, we need the following maximum principles, which are known for experts (see for examples \cite{GLY16,GLY16'}). For readers' convenience, we include the detailed proofs here.
\begin{lemma}[Weak maximum principle]\label{lemma:weak}
For any constant $c>0$, if $u$ satisfies $\Delta u+cu\geq 0$, then $u\geq 0$ on $V$.
\end{lemma}
\begin{proof}
Let $u^-=\min\{u,0\}$. For any $x\in V$, we \textit{claim }that
\begin{equation}\label{eq:maxi}
\Delta u^-(x)+cu^-(x)\geq0,
\end{equation}
from which, one has
\[\int_V\Gamma(u^-)d\mu+c\|u^-\|^2_{\ell^2_\mu}=\langle u^-,\Delta u^-+cu^-\rangle\leq 0.\]
This leads to $u^-\equiv0$ on $V$.

To prove this claim, we first consider the case $u(x)\geq 0$. Therefore, $cu^-(x)=0$ and
\[\Delta u^-(x)=\frac{1}{\mu(x)}\sum_{y\sim x}\omega_{xy}(u^-(x)-u^-(y))=-\frac{1}{\mu(x)}\sum_{y\sim x}\omega_{xy}u^-(y)\geq 0,\]
due to $u^-(z)\leq0$ for any $z\in V$.
In the case $u(x)<0$, one has $cu^-(x)=cu(x)$ and thus
\[\begin{split}
\Delta u^-(x)=\frac{1}{\mu(x)}\sum_{y\sim x}\omega_{xy}(u^-(x)-u^-(y))
&=\frac{1}{\mu(x)}\sum_{y\sim x}\omega_{xy}(u(x)-u^-(y))\\
&\geq\frac{1}{\mu(x)}\sum_{y\sim x}\omega_{xy}(u(x)-u(y))=\Delta u(x).
\end{split}\]
It follows that $\Delta u^-(x)+cu^-(x)\geq \Delta u(x)+cu(x)\geq 0$, which confirms \eqref{eq:maxi}
and ends the proof of the lemma.
\end{proof}
\begin{lemma}[Strong maximum principle]\label{lemma:strong}
Suppose that $u\geq 0$, and that $\Delta u+cu\geq 0$ for some constant $c>0$.
If there exists $x_0\in V$ such that $u(x_0)=0$, then $u\equiv0$ on $V.$
\end{lemma}
\begin{proof}
Let $x=x_0$, we have
\[\frac{1}{\mu(x_0)}\sum_{y\sim x_0}\omega_{yx_0}(u(x_0)-u(y))+cu(x_0)\geq0,\]
which implies
\[\frac{1}{\mu(x_0)}\sum_{y\sim x_0}\omega_{yx_0}u(y)\leq0.\]
Since $u\geq 0$ and $\omega_{yx_0}>0$ for all $y\sim x_0$, we obtain
\[u(y)=0,~~~~\mbox{for all}~~y\sim x_0.\]
Therefore, $u\equiv0$ on $V$ by the connectedness of $G.$
\end{proof}
\subsection{Palais-Smale condition}
We define a functional $E_\lambda:W^{1,2}(V)\rightarrow \mathbb {R}$ by
\begin{equation*}\label{eq:va}
E_\lambda(u)=\int_V(|\nabla u|^2+2\kappa u-K_\lambda e^{2u})d\mu,
\end{equation*}
where $\kappa$ and $K_\lambda$ are given as in the assumptions of Theorem \ref{thm:main}, in particular
$\int_V\kappa d\mu<0$. For any $\phi\in W^{1,2}(V)$, denote by $dE_\lambda(u)(\phi)$ the Frechet derivative
of the functional, by $d^kE_\lambda(u)(\phi,\cdots,\phi)$ the Frechet derivative of order $k\geq 2$.
\begin{lemma}[Palais-Smale condition]\label{lemma:ps}
Suppose that $V_\lambda^-=\{x\in V:K_\lambda(x)<0\}$ is nonempty for some $\lambda\in \mathbb{R}$. Then $E_\lambda$ satisfies the $(PS)_c$ condition for all $c\in \mathbb{R}$, i.e. if $(u_j)$ is a sequence of functions in $W^{1,2}(V)$ such that $E_\lambda(u_j)\rightarrow c$ and $dE_\lambda(u_j)(\phi)\rightarrow 0$ for all $\phi\in W^{1,2}(V)$ as $j\rightarrow\infty$, then there exists some $u_0\in W^{1,2}(V)$ satisfying $u_j \rightarrow u_0$ in $W^{1,2}(V).$
\end{lemma}
\begin{proof}
Let $(u_j)$ be a function sequence such that $E_\lambda(u_j)\rightarrow c$ and $dE_\lambda(u_j)(\phi)\rightarrow 0,$ or equivalently
\begin{equation}\label{eq:1}
\int_V(|\nabla u_j|^2+2\kappa u_j-K_\lambda e^{2u_j})d\mu=c+o_j(1),
\end{equation}
\begin{equation}\label{eq:2}
\int_V(\nabla u_j\nabla \phi+\kappa \phi-K_\lambda e^{2u_j}\phi)d\mu=o_j(1)\|\phi\|_{W^{1,2}(V)}, ~~~~\forall \phi\in W^{1,2}(V),
\end{equation}
where $o_j(1)\rightarrow 0$ as $j\rightarrow \infty.$

Let $\phi\equiv1$ in \eqref{eq:2}, one has
\[\int_V(\kappa-K_\lambda e^{2u_j})d\mu=o_j(1)\mu(V)^{1/2},\]
which implies
\begin{equation}\label{eq:3}
\int_VK_\lambda e^{2u_j}d\mu=\int_V\kappa d\mu+o_j(1).
\end{equation}
Inserting \eqref{eq:3} into \eqref{eq:1}, we obtain
\begin{equation}\label{eq:4}
\int_V(|\nabla u_j|^2+2\kappa u_j)d\mu=\int_V \kappa d\mu+c+o_j(1).
\end{equation}

We now \textit{claim} that $u_j$ is bounded in $\ell^2_\mu$. Suppose not, there holds $\|u_j\|_{\ell^2_\mu}\rightarrow \infty$. We set $v_j=\frac{u_j}{\|u_j\|_{\ell^2_\mu}}.$  By the Cauchy-Schwarz inequality, one has
\[\int_V\kappa\frac{u_j}{\|u_j\|^2_{\ell^2_\mu}}d\mu=o_j(1).\]
This together with \eqref{eq:4} leads to
\begin{equation}\label{eq:5}
\int_V|\nabla v_j|^2d\mu=o_j(1).
\end{equation}
Hence, $v_j$ is bounded in $W^{1,2}(V)$. In view of Lemma \ref{lemma:compact} and \eqref{eq:5}, $v_j\rightarrow \gamma$ in $W^{1,2}(V)$ for some constant $\gamma.$ Here and in the sequel, we do not distinguish sequence and subsequence. Since $\|v_j\|_{\ell^2_\mu}=1$, we have $\gamma\neq0.$ It follows from \eqref{eq:4} that
\[\int_V\kappa v_jd\mu\leq o_j(1).\]
Passing to the limit $j\rightarrow\infty$ in the above inequality, we conclude that $\gamma\geq 0$ since $\int_V\kappa d\mu< 0$. Therefore $\gamma>0$.

On the other hand, for any $x\in V_\lambda^-$, if there exists $N\in \mathbb{N}$, if $j>N$ such that $u_{j}(x)\leq 0$, then $\lim_{j\rightarrow \infty}v_j(x)\leq 0$, which contradicts $\gamma>0$ and confirms our claim. If not, let $x_*\in V_\lambda^-$, due to the finiteness of $V$, we can choose a subsequence $\{j_k\}_{k=0}^\infty$ such that $u_{j_k}(x_*)>0$. Set
\[\phi(x)=\left\{
            \begin{array}{ll}
              u_{j_k}(x_*), & \hbox{$x=x_*$} \\[1.2ex]
              0, & x\not=x_\ast.
            \end{array}
          \right.
\]
Then
\[\|\phi\|_{W^{1,2}(V)}^2
=2\sum_{y\sim x_*}\omega_{x_*y}u_{j_k}^2(x_*)+\mu(x_*)u_{j_k}^2(x_*)
=(2\deg(x_*)+\mu(x_*))u_{j_k}^2(x_*).
\]
Substituting it into \eqref{eq:2}, we have
\begin{equation}\label{eq:con}
\Delta u_{j_k}(x_*)+\kappa(x_*)-K_\lambda(x_*) e^{2u_{j_k}(x_*)}\leq C'.
\end{equation}
Since $v_j\ra\gamma$, $\|u_j\|_{\ell^2_\mu}\ra+\infty$ and $V$ has finite points, we conclude
$$u_j=(\gamma+o_j(1))\|u_j\|_{\ell^2_\mu}\quad{\rm uniformly\,\,on}\quad V.$$
This together with (\ref{eq:con}) leads to
\[\|u_{j_k}\|_{\ell^2_\mu}^2o_{j_k}(1)+\kappa(x_*)-K_\lambda(x_*) e^{2(\gamma+o_{j_k}(1))\|u_{j_k}\|_{\ell^2_\mu}^2}\leq C',\]
which is impossible since $K_\lambda(x_\ast)<0$. Then our claim follows immediately.

Since $u_j$ is bounded in $\ell_\mu^2$, we have $u_j$ is bounded in $W^{1,2}(V)$ due to the finiteness of $V$. Therefore, by Lemma \ref{lemma:compact}, there exists some $u_0\in W^{1,2}(V)$ such that up to subsequence, $u_j\rightarrow u_0$ in $W^{1,2}(V)$.
\end{proof}
\subsection{Upper and lower solutions principle}
Let $f:V\times \mathbb{R}\rightarrow \mathbb{R}$ be a function, and $f$ is smooth with respect to the second variable. We say that $u\in V^\mathbb{R}$ is an upper (lower) solution to the following equation
\begin{equation}\label{eq:elliptic}
\Delta u(x)+f(x,u(x))=0,~~x\in V,
\end{equation}
if $u$ satisfies $\Delta u(x)+f(x,u(x))\geq (\leq) ~0$ for any $x\in V$. We generalize (\cite{GLY16}, Lemma 8) to
the following:
\begin{lemma}\label{lemma:upper}
Suppose that $\varphi,\psi$ are lower and upper solution to \eqref{eq:elliptic} respectively with $\varphi\leq\psi$ on $V$. Then \eqref{eq:elliptic} has a solution $u$ with $\varphi\leq u\leq\psi$ on $V$.
\end{lemma}

\begin{proof} This is a discrete version of the argument of Kazdan-Warner (\cite{KW74}, Lemma 9.3), and the method of proof carries over to the setting of graphs.

 Since the graph is finite, there exists a constant $A$ such that $-A\leq \varphi\leq \psi\leq A$. One can find a sufficient large constant $c$ such that $F(x,t)=ct-f(x,t)$ is increasing with respect to $t\in[-A,A]$ for any fixed $x\in V$. We define an operator $Lu=\Delta u+cu$, and $L$ is a compact operator and $\mbox{Ker}(L) = \mbox{span}\{1\}$ due to the finiteness of the graph. Hence, we can define $\varphi_{j+1},\psi_{j+1}$ inductively as the unique solution to
 \[\varphi_0=\varphi, ~~L\varphi_{j+1}(x)=c\varphi_j(x)-f(x,\varphi_j(x)), ~~\forall j\geq0, x\in V,\]
 \[\psi_0=\psi, ~~L\psi_{j+1}(x)=c\psi_j(x)-f(x,\psi_j(x)), ~~\forall j\geq0, x\in V\]
 respectively.
Combining with the definition of upper (lower) solution and the monotonicity of $F(x,t)$ with respect to $t$, we obtain
\[L\varphi_0(x)\leq L\varphi_1(x)=F(x,\varphi(x))\leq F(x,\psi(x))=L\psi_1(x)\leq L\psi(x),~~ x\in V.\]
Then the weak maximum principle (see Lemma \ref{lemma:weak}) yields that
\[\varphi\leq \varphi_1\leq\psi_1\leq \psi.\]
Moreover, it turns out that $\varphi_1$ and $\psi_1$ are lower and upper solution to \eqref{eq:elliptic} respectively. By induction, we have
\[\varphi\leq \varphi_{j}\leq \varphi_{j+1}\leq\psi_{j+1}\leq\psi_j\leq \psi, ~~j=1,2,\cdots.\]
Since $V$ is finite, it is
easy to see that up to a subsequence, $\varphi_j \rightarrow u_1,\psi_j\rightarrow u_2$ uniformly on $V$, and $u=u_1$ or $u_2$ is a solution to \eqref{eq:elliptic} with $\varphi\leq u\leq \psi$ on $V$.
\end{proof}

\section{Proof of Theorem \ref{thm:main}}\label{sec-3}
\subsection {Unique solution in the case $\lambda\leq0$.}
\begin{claim}\label{claim1}
$E_\lambda$ is strictly convex on $W^{1,2}(V)$.
\end{claim}
\begin{proof}
We only need to show that there exists some constant $C>0$ such that
\begin{equation}\label{eq:sconvex}
d^2 E_\lambda(u)(h,h)\geq C\|h\|^2_{W^{1,2}(V)},~~\forall u, h\in W^{1,2}(V).
\end{equation}
Suppose not, there would be a function $u\in W^{1,2}(V)$ and a function sequence $h_j\in W^{1,2}(V)$ such that $\|h_j\|_{W^{1,2}(V)}=1$ for all $j$ and $d^2 E_\lambda(u)(h_j,h_j)\rightarrow 0$ as $j\rightarrow \infty$. From Lemma \ref{lemma:compact}, there exists $h_\infty\in W^{1,2}(V)$, such that up to a subsequence, $h_j\rightarrow h_\infty$ as $j\rightarrow \infty$ in $W^{1,2}(V)$. Since
\[d^2E_\lambda(u)(h_j,h_j)=2\int_V(|\nabla h_j|^2-2K_\lambda e^{2u}h_j^2)d\mu,\]
and $K_\lambda\leq 0$, it follows that $\int_V|\nabla h_j|^2d\mu\rightarrow0$ and $\int_VK_\lambda e^{2u}h_j^2d\mu\rightarrow0$, which lead to $h_\infty\equiv c$ for some constant $c$, and moreover
\[c^2\int_VK_\lambda e^{2u}d\mu=\lim_{j\rightarrow \infty}\int_VK_\lambda e^{2u}h_j^2d\mu=0.\]
It is easily seen that $\int_VK_\lambda e^{2u}d\mu<0$ by $K\not\equiv0$, thus $c=0$. This contradicts \[\|h_\infty\|_{W^{1,2}(V)}=\lim_{j\rightarrow \infty}\|h_j\|_{W^{1,2}(V)}=1.\]
 Hence \eqref{eq:sconvex} holds.
\end{proof}
\begin{claim}\label{claim2}
For any $\varepsilon>0$, there exist constants $C, C(\varepsilon)>0$ such that
\[E_\lambda(u)\geq (C-2\varepsilon)\|u\|_{W^{1,2}(V)}^2-2C(\varepsilon).\]
\end{claim}
\begin{proof}
By Young's inequality, for any $\varepsilon>0$, there exists a constant $C(\varepsilon)>0$ such that
\[\le|\int_V\kappa ud\mu\ri|\leq \varepsilon\|u\|_{W^{1,2}(V)}^2+C(\varepsilon).\]
Thus, it is sufficient to find some constant $C>0$ such that for all $u\in W^{1,2}(V)$
\[\int_V(|\nabla u|^2-K_\lambda e^{2u})d\mu\geq C\|u\|^2_{W^{1,2}(V)}.\]
Suppose not, there would exist a sequence of functions $u_j$ satisfying
\[\int_V(|\nabla u_j|^2+u_j^2)d\mu=1,~~\int_V(|\nabla u_j|^2-K_\lambda e^{2u_j})d\mu=o_j(1).\]
Clearly, $u_j$ is bounded in $W^{1,2}(V)$, it follows from Lemma \ref{lemma:compact} that there exists some function $u\in W^{1,2}(V)$ such that up to a subsequence, $u_j\rightarrow u$ in $W^{1,2}(V)$ as $j\rightarrow \infty.$ Due to $K_\lambda\not\equiv0$ and $K_\lambda\leq 0,$ we have
\[0<\int_V(|\nabla u|^2-K_\lambda e^{2u})d\mu=\lim_{j\rightarrow \infty}\int_V(|\nabla u_j|^2-K_\lambda e^{2u_j})d\mu=0,\]
which gets a contradiction.
\end{proof}
\textit{Proof of (1) in Theorem \ref{thm:main}.}
~~It is a consequence of Claim \ref{claim1} and Claim \ref{claim2}. Precisely we denote $\Lambda=\inf_{u\in W^{1,2}(V)} E_\lambda(u)$. By Claim \ref{claim2},
we see that $\Lambda$ is a definite real number. Take a function sequence $u_j\in W^{1,2}(V)$ such that $E_\lambda(u_j)\ra \Lambda$ as $j\ra\infty$. Applying Claim \ref{claim2}, we have that $u_j$ is bounded in $W^{1,2}(V)$. Then, in view of Lemma \ref{lemma:compact},
there exists a subsequence of $u_j$ (still denoted by $u_j$) and a function $u_0\in W^{1,2}(V)$ such that $u_j\ra u_0$ in $W^{1,2}(V)$. Obviously
$E_\lambda(u_0)=\Lambda$, and thus $u_0$ is a critical point of $E_\lambda$. We also need to explain why $E_\lambda$ has only one critical point. For otherwise, we assume $u^\ast$ is another critical point of $E_\lambda$. Note that $dE_\lambda(u_0)=dE_\lambda(u^\ast)=0$. It follows from
Claim \ref{claim1}, particularly from (\ref{eq:sconvex}), that
\bna
E_\lambda(u_0)&=&E_\lambda(u^\ast)+dE_\lambda(u^\ast)(u_0-u^\ast)+\f{1}{2}d^2E_\lambda(\xi)(u_0-u^\ast,u_0-u^\ast)\\
&\geq& E_\lambda(u^\ast)+C\|u_0-u^\ast\|_{W^{1,2}(V)}^2\ena
for some positive constant $C$,
where $\xi$ is a function lies between $u^\ast$ and $u_0$.
Hence we have $E_\lambda(u_0)>E_\lambda(u^\ast)$, contradicting the fact that $E_\lambda(u_0)=\Lambda$. This implies the uniqueness of the critical
point of $E_\lambda$. $\hfill\Box$

\subsection{Multiplicity of solutions for $0<\lambda<\lambda^*$.}
 Fixing $\lambda\in(0,\lambda^*)$, we will seek two different solutions of \eqref{eq:0}. One is a strict local minimum of the functional $E_\lambda$, and the other is from mountain-pass theorem.
We firstly prove the existence of $\lambda^*$. Consider the case $\lambda=0$ in the equation \eqref{eq:0} as follows
\begin{equation}\label{eq:00}
\Delta u+\kappa-Ke^{2u}=0.
\end{equation}
For the solution $u_0$ of  \eqref{eq:00},  the linearized equation of \eqref{eq:00} at $u_0$
\[\Delta v-2Ke^{2u_0}v=0\]
has only a trivial solution $v\equiv0,$ since $K\leq 0$ and $K\not\equiv0.$ Indeed,
\[0\leq \int_V|\nabla v|^2d\mu=\int_V v\Delta vd\mu=2\int_VKe^{2u_0}v^2d\mu\leq0,\]
which implies $v(x)=0$ when $K(x)<0.$ In the case $K(x)=0$, we have $\Delta v(x)=0$. Thus
\[\int_V|\nabla v|^2d\mu=\int_V v\Delta vd\mu=0.\]
It follows that $v$ is a constant function and hence $v\equiv0.$
By the implicit function theorem, there exists a small enough $s>0$ such that the equation \eqref{eq:0} has a solution for any $\lambda\in(0,s)$. Indeed, let $u=tv+u_0, v\not\equiv 0$ on $V$, we consider $G(\lambda,t)=\Delta u+\kappa-K_\lambda e^u$. It is easy to see that $G(0,0)=0$, $G(\lambda,t)$ and $\partial_tG(\lambda,t)=\Delta v-2v(K+\lambda) e^{2(u_0+tv)}$ are continuous on any domain $D\subset\mathbb{R}^2$, Furthermore,  $\partial_tG(0,0)=\Delta v-2K e^{2u_0}v\not\equiv 0$ unless $v\equiv 0$ on $V$. Therefore, by the implicit function theorem, there exists $s>0$ such that $t=g(\lambda)$ and $G(\lambda,g(\lambda))=0$ for any $\lambda\in(0,s)$.  In other words, $u_\lambda=g(\lambda)v+u_0, \forall \lambda\in(0,s)$ is the solution of  \eqref{eq:0}. Define
\[\lambda^*=\sup\{s: \mbox{the equation \eqref{eq:0} has a solution for any $\lambda\in (0,s)$}\}.\]
One can see that $\lambda^*\leq -\min_V K.$ For otherwise, $K_\lambda=K+\lambda>0$ for some $\lambda<\lambda^*$. Adding up the equation \eqref{eq:0} for all $x\in V$, we have
\[0>\int_V\kappa d\mu=\int_VK_\lambda e^{2u}d\mu>0,\]
which is impossible. In conclusion, we have $0<\lambda^*\leq -\min_V K$.

\textit{Proof of (2) in Theorem \ref{thm:main}.} We separate the proof into the
following three steps.

\textbf{Step 1.} \textit{The existence of the upper and lower solution of \eqref{eq:0}.}

Take $\lambda_1$ with $\lambda<\lambda_1<\lambda^*$, let $u_{\lambda_1}$ be a solution of  \eqref{eq:0} at $\lambda_1$. It is easily seen that $\psi=u_{\lambda_1}$ is a strict upper solution of \eqref{eq:0} at $\lambda$, namely
\[\Delta \psi+\kappa-K_\lambda e^{2\psi}>0.\]
Let $v$ be the solution to the following equation
\begin{equation}\label{eq:lo-solu}
\Delta v=-\kappa+\frac{1}{\mu(V)}\int_V\kappa d\mu.
\end{equation}
The existence of solution to \eqref{eq:lo-solu} was proved in \cite{GLY16}. Set $\varphi=v-s$, where $s$ is a sufficiently large constant such that $\varphi<\psi$ on $V$ and
\[\Delta \varphi+\kappa-K_\lambda e^{2\varphi}=\frac{1}{\mu(V)}\int_V\kappa d\mu-K_\lambda e^{2v-2s}<0\]
since $\int_V\kappa d\mu<0$. Therefore, $\varphi$ is a strict lower solution of \eqref{eq:0}. Let $[\varphi,\psi]$ be the order interval defined by
\[[\varphi,\psi]=\{u\in V^\mathbb{R}:\varphi\leq u\leq\psi~\mbox{on}~ V\}.\]
The upper and lower-solution method (Lemma \ref{lemma:upper}) asserts that \eqref{eq:0} has a solution $u_\lambda\in[\varphi,\psi]$ on $V$.

\textbf{Step 2.} \textit{$u_\lambda$ can be chosen as a strict local minimum of $E_\lambda$.}

Let $f_\lambda(x,t)=ct-\kappa(x)+K_\lambda(x) e^{2t}$, where $c$ is sufficiently large such that $f_\lambda(x,t)$ is increasing in $t\in[-A,A]$, $A$ is a constant such that $-A\leq \varphi<\psi\leq A$ on $V$. Let $F_\lambda(x,u(x))=\int_0^{u(x)}f_\lambda(x,t)dt$. It is easy to rewrite $E_\lambda(u)$ as
\[E_\lambda (u)=\int_V(|\nabla u|^2+cu^2)d\mu-2\sum_{x\in V}\mu(x)F_\lambda(x,u(x))-\int_VK_\lambda d\mu.\]

It is obvious that $E_\lambda$ is bounded from below on $[\varphi,\psi]$. 
Therefore, we denote
\[a:=\inf_{u\in [\varphi,\psi]}E_\lambda(u).\]
Taking a function sequence $u_j\subset [\varphi,\psi]$ such that $ E_\lambda(u_j)\rightarrow a$ as $j\rightarrow \infty.$ From it, we can get that $u_j$ is bounded in $W^{1,2}(V)$, and thus up to subsequence, $u_j$ converges to some $u_\lambda$ in $W^{1,2}(V)$ and $\ell^q_\mu$ for any $q\geq 1$, and $e^{2u_j}$ converges to $e^{2u_\lambda}$ in $\ell^1_\mu$. Hence
\[E_\lambda(u_\lambda)=\inf_{u\in [\varphi,\psi]} E_\lambda(u).\]
As a consequence, $u_\lambda$ satisfies the Euler-Lagrange equation
\[\Delta u_\lambda(x)+cu_\lambda(x)=f_\lambda(x,u_\lambda(x)).\]
From it, one can conclude that
\begin{equation}\label{eq:esti}
\varphi(x)<u_\lambda(x)<\psi(x), \forall x\in V.
\end{equation}
Indeed, noting that $f_\lambda(x,t)$ is increasing with respect to $t\in[-A,A]$, we have
\[\Delta \varphi(x)+c\varphi(x)\leq f_\lambda(x,\varphi(x))\leq f_\lambda(x,\psi(x))\leq \Delta \psi(x)+c\psi(x).\]
One can conclude \eqref{eq:esti} by the strong maximum principle (Lemma \ref{lemma:strong}), and the fact $\varphi<\psi$.
For any $h\in W^{1,2}(V)$, we define a function $\eta(t)=E_\lambda(u_\lambda+th), t\in\mathbb{R} $.  There holds $\varphi\leq u_\lambda+th\leq \psi$ for sufficiently small $|t|$.
Since $u_\lambda$ is a minimum of $E_\lambda$ on $(\varphi,\psi)$, we have $\eta'(0)=dE_\lambda(u_\lambda)(h)=0$ and $\eta''(0)=d^2E_\lambda(u_\lambda)(h,h)\geq0.$ 
Furthermore, there exists a constant $C>0$ such that
\begin{equation}\label{eq:lower}
d^2E_\lambda(u_\lambda)(h,h)\geq C\|h\|_{W^{1,2}(V)}, ~~~~\forall h\in W^{1,2}(V),
\end{equation}
which implies  $u_\lambda$ is strict local minimum of $E_\lambda$ on $W^{1,2}(V)$.
It remains to prove \eqref{eq:lower}. We first denote
\[\theta:=\inf_{\|h\|_{W^{1,2}(V)}=1}d^2E_\lambda(u_\lambda)(h,h),\]
which is nonnegative. It is sufficient to prove $\theta>0$, \eqref{eq:lower} follows.
Suppose $\theta=0$, we \textit{claim} that there exists some $\tilde{h}$ with $\|\tilde{h}\|_{W^{1,2}(V)}=1$ such that $d^2E_\lambda(u_\lambda)(\tilde{h},\tilde{h})=0$. To see this, let $h_j$ be a function sequence satisfying  $\|h_j\|_{W^{1,2}(V)}=1$ for all $j$ and $d^2E_\lambda(u_\lambda)(h_j,h_j)\rightarrow 0$ as $j\rightarrow \infty.$  Up to subsequence, $h_j\rightarrow \tilde{h}$ in $W^{1,2}(V)$ from Lemma \ref{lemma:compact}, and confirms our claim. To put it another way, the functional $v\mapsto d^2E_\lambda(u_\lambda)(v,v)$ attains its minimum at $v=\tilde{h}$, it follows that $d^2E_\lambda(u_\lambda)(\tilde{h},v)=0$
for all $v\in W^{1,2}(V).$
 Hence, $\tilde{h}$ is a solution of the following equation
\begin{equation}\label{eq:eqution}
\Delta h=2K_\lambda e^{2u_\lambda}h,~~~~\lambda\in(0,\lambda^*).
\end{equation}
It is easy to see that $\tilde{h}$ is not a constant. For otherwise \eqref{eq:eqution} yields
\[0>\int_V \kappa d\mu=\int_VK_\lambda e^{2u_\lambda}d\mu=0,\]
which is impossible. Multiplying \eqref{eq:eqution} by $\tilde{h}^3$, we obtain
\[\begin{split}
d^4E_\lambda(u_\lambda)(\tilde{h},\tilde{h},\tilde{h},\tilde{h})
&=-16\int_VK_\lambda e^{2u_\lambda}\tilde{h}^4d\mu\\
&=-8\int_V\tilde{h}^3\Delta \tilde{h}d\mu\\
&=-8\sum_{x,y\in V}\omega_{xy}(\tilde{h}^3(x)-\tilde{h}^3(y))(\tilde{h}(x)-\tilde{h}(y))\\
&=-8\sum_{x,y\in V}\omega_{xy}(\tilde{h}^2(x)+\tilde{h}(x)\tilde{h}(y)+\tilde{h}^2(y))(\tilde{h}(x)-\tilde{h}(y))^2<0.
\end{split}\]
The last inequality is due to the fact $a^2+ab+b^2>0$ for any $a,b\in \mathbb{R}$ satisfying $ab\neq0$.
Since $d^2E_\lambda(u_\lambda+t\tilde{h})(\tilde{h},\tilde{h})$ attains its minimum at $t=0$, we have $d^3E_\lambda(u_\lambda)(\tilde{h},\tilde{h},\tilde{h})=0$, which together with $dE_\lambda(u_\lambda)(\tilde{h})=0$ and $d^2E_\lambda(u_\lambda)(\tilde{h},\tilde{h})=0$ leads to
\begin{equation}\label{eq:estim}
E_\lambda(u_\lambda+\epsilon \tilde{h})=E_\lambda(u_\lambda)+\frac{\epsilon^4}{24}d^4E_\lambda(u_\lambda)(\tilde{h},\tilde{h},\tilde{h},\tilde{h})+0(\epsilon^5)<E_\lambda(u_\lambda)
\end{equation}
for small $\epsilon>0.$ Let $\epsilon$ small enough such that $\varphi\leq u_\lambda+\epsilon \tilde{h}\leq \psi$, thus by \eqref{eq:estim},
\[E_\lambda(u_\lambda+\epsilon \tilde{h})<E_\lambda(u_\lambda),\]
which contradicts the fact that $u_\lambda$ is the minimum of $E_\lambda$ on $[\varphi,\psi]$. Therefore $\theta>0$, which concludes \eqref{eq:lower}.

\textbf{Step 3.} \textit{ The second solution of \eqref{eq:0} is given by the mountain-pass theorem.}

 We shall use the mountain-pass theorem due to Ambrosetti and Rabinowitz \cite{AR73}, which reads as follows:
 {\it Let
 $(X,\|\cdot\|)$ be a Banach space, $J\in C^1(X,\mathbb{R})$, $e_0,e\in X$ and
 $r>0$ be such that $\|e-e_0\|>r$ and
 $$b:=\inf_{\|u-e_0\|=r}J(u)>J(e_0)\geq J(e).$$
 If $J$ satisfies the $(PS)_c$ condition with $c:=\inf_{\gamma\in\Gamma}\max_{t\in[0,1]}J(\gamma(t))$,
 where $$\Gamma:=\{\gamma\in C([0,1],X): \gamma(0)=e_0,\gamma(1)=e\},$$
 then $c$ is a critical value of $J$.} In our case, $W^{1,2}(V)$ is a Banach space, and $E_\lambda:W^{1,2}(V)\ra\mathbb{R}$ is a
 smooth functional.

 Since $u_\lambda$ is a strict local minimum of $E_\lambda$ on $W^{1,2}(V)$, there exists a small enough number $r>0$ such that
\begin{equation}\label{eq:m1}
\inf_{\|u-u_\lambda\|_{W^{1,2}(V)}=r}E_\lambda(u)>E_\lambda(u_\lambda).
\end{equation}
Moreover, for any $\lambda>0,$ $E_\lambda$ has no lower bound on $W^{1,2}(V)$, namely, there exists $v\in W^{1,2}(V)$ such that
\begin{equation}\label{eq:m2}
E_\lambda(v)<E_\lambda(u_\lambda), ~~~~ \|v-u_\lambda\|_{W^{1,2}(V)}>r.
 \end{equation}
To see this, we set $V_\varepsilon=\{x\in V: K_\lambda(x)>\varepsilon\}$ for small $\varepsilon>0$. Note that $V_\varepsilon$ is nonempty since $\max_{ V}K=0$ and $\lambda>0$. Let $f\in W^{1,2}(V)$ be a function which equals to $1$ in $V_\varepsilon$ and vanishes on  $V/V_\varepsilon$, then
\[\begin{split}
E_\lambda(tf )&=t^2\int_V|\nabla f|^2d\mu+t\int_{V_\varepsilon}\kappa d\mu-\int_{V_\varepsilon}K_\lambda e^{2t}d\mu-\int_{V/V_\varepsilon}K_\lambda d\mu\\
&\leq At^2+Bt+C-\varepsilon \mu(V_\varepsilon) e^{2t}\rightarrow -\infty,~~~~ t\rightarrow +\infty.
\end{split}\]
In view of \eqref{eq:m1}, \eqref{eq:m2} and Lemma \ref{lemma:ps}, the mountain pass theorem of Ambrosetti and Rabinowitz gives another critical point $u^\lambda$ of $E_\lambda$ other than $u_\lambda$. In particular,
\[E_\lambda(u^\lambda)=\min_{\gamma\in \Gamma}\max_{u\in \gamma}E_\lambda(u),\]
where $\Gamma=\{\gamma\in C([0,1],W^{1,2}(V)):\gamma(0)=u_\lambda,\gamma(1)=v\}$, and $dE_\lambda(u^\lambda)(h)=0$ for any $h\in W^{1,2}(V)$.
$\hfill\Box$

\subsection{Solvability at $\lambda^*$.}

For any $\lambda, 0<\lambda<\lambda^*$, let $u_\lambda$ be the local minimum of $E_\lambda$ obtained in the previous subsection. That is, $u_\lambda$ is the solution of \eqref{eq:0}, and
\begin{equation}\label{eq:convex}
d^2E_\lambda(u_\lambda)(h,h)=2\int_V(|\nabla h|^2-2K_\lambda e^{2u_\lambda}h^2)d\mu\geq0,~~~~\forall h\in W^{1,2}(V).
\end{equation}
\textit{Proof of (3) in Theorem \ref{thm:main}.}
The crucial point in this proof is to show that $u_\lambda$ is uniformly bounded in $W^{1,2}(V)$ as $\lambda\rightarrow \lambda^*.$ If it is true, then up to subsequence, $u_\lambda$ converges to some $u$ in $W^{1,2}(V)$, and $u$  is the solution of
\begin{equation*}\label{lambda}
\Delta u+\kappa-K_{\lambda^*}e^{2u}=0.
\end{equation*}
Hence, we aim to prove the $W^{1,2}(V)$ boundedness of $u_\lambda.$ To this end, we divide the proof into the following three steps.

\textbf{Step 1.} \textit{There exists a constant $C>0$ such that
\begin{equation}\label{lower}
u_\lambda\geq-C~~~~\mbox{on $V$}
\end{equation}
uniformly for any $0<\lambda<\lambda^*$.}

Let $v$ satisfy \eqref{eq:lo-solu}, and $\varphi_s=v-s$ for $s>0$. Then for sufficient large $s$, say $s\geq s_0$, $\varphi_s$ is a continuous family with respect to $s$ of strict lower solution of \eqref{eq:0} at $\lambda=0$, i.e.
\[\Delta \varphi_s+\kappa-K e^{2\varphi_s}<0.\]
It is clear that $\varphi_s$ is also a strict lower solution of \eqref{eq:0} at $\lambda\in (0,\lambda^*)$ for any $s\in [s_0,\infty)$.
Now, we prove that $u_\lambda\geq\varphi_{s_0}$, and thus \eqref{lower} holds. For otherwise, we can find for some $s\in (s_0,\infty)$
\[u_\lambda\geq \varphi_s~~~\mbox{on $V$}, ~~~~\mbox{and}~~u_\lambda(\tilde{x})=\varphi_s(\tilde{x})~~\mbox{for some $\tilde{x}\in V$}.\]
It follows from the strong maximum principle (Lemma \ref{lemma:strong}) that $u_\lambda\equiv\varphi_s$, which contradicts that $\varphi_s$ is strict low solution of \eqref{eq:0} at $\lambda\in[0,\lambda^*).$

\textbf{Step 2.} \textit{The set $V_{\lambda^*}^-=\{x\in V: K_{\lambda^*}(x)<0\}$ is not empty. }

From the case (1) in Theorem \ref{thm:main}, there exists unique solution $w_0$ of the equation
\[\Delta w+\kappa+e^{2w}=0.\]
Together with the solution $u_\lambda$ to the equation \eqref{eq:0} at $\lambda$, and let $v_\lambda=u_\lambda-w_0$, we have
\[\Delta v_\lambda-K_\lambda e^{2u_\lambda}-e^{2w_0}=0.\]
Multiplying the above equation by $e^{-2v_\lambda}$ and intergrading by parts, one has
\[\int_V K_\lambda e^{w_0}d\mu=\int_Ve^{-2v_\lambda}\Delta v_\lambda d\mu-\int_Ve^{-2(u_\lambda-2w_0)}d\mu\leq0.\]
Indeed, by Green formula
\[\int_Ve^{-2v_\lambda}\Delta v_\lambda d\mu=\sum_{x,y\in V}\omega_{xy}(v_\lambda(x)-v_\lambda(y))(e^{-2v_\lambda(x)}-e^{-2v_\lambda(y)})\leq0\]
since $(v_\lambda(x)-v_\lambda(y))(e^{-2v_\lambda(x)}-e^{-2v_\lambda(y)})\leq0$ for any $x,y\in V$. Therefore
\[\int_V K_{\lambda^*}e^{w_0}d\mu=\lim_{\lambda\rightarrow\lambda^*}\int_{V}K_\lambda e^{w_0}d\mu\leq 0.\]
Suppose that $K_{\lambda^*}\geq 0,$ thus $K_{\lambda^*}\equiv 0$. This contradicts the assumption that $K_{\lambda^*}$ is not a constant.

\textbf{Step 3.} \textit{$u_\lambda$ is uniformly bounded in $W^{1,2}(V)$ as $\lambda\rightarrow \lambda^*$.}

Fixing $x_0\in V_{\lambda^*}^-$, we set $\rho\in V^{\mathbb{R}}$ to be a function which vanishes besides $x_0$ and $\rho(x_0)<0$. Consider the equation
\[\Delta w+\kappa-\rho e^{2w}=0,\]
which always has the unique solution from the case (1) in Theorem \ref{thm:main}, set $w_1.$ As above, the function $v_\lambda=u_\lambda-w_1$ satisfies the equation
\begin{equation}\label{eq:eq}
\Delta v_\lambda+\rho e^{2w_1}-K_\lambda e^{2(v_\lambda+w_1)}=0.
\end{equation}
Multiplying \eqref{eq:eq} by $e^{2v_\lambda}$ and integrating by parts over $V$ gives
\begin{equation}\label{eq:eq1}
\int_Ve^{2v_\lambda}\Delta v_\lambda d\mu+\int_V\rho e^{2(v_\lambda+w_1)}d\mu-\int_VK_\lambda e^{2(2v_\lambda+w_1)}d\mu=0.
\end{equation}
Utilizing the fact that
\[(e^{2a}-e^{2b})(a-b)\geq (e^a-e^b)^2,~~~~a,b\in \mathbb{R},\]
one can estimate the first term in \eqref{eq:eq1} as follows.
\[\begin{split}
\int_Ve^{2v_\lambda}\Delta v_\lambda d\mu
&=\frac{1}{2}\sum_{x,y\in V}\omega_{xy}(e^{2v_\lambda(x)}-e^{2v_\lambda(y)})(v_\lambda(x)-v_\lambda(y))\\
&\geq \frac{1}{2}\sum_{x,y\in V}\omega_{xy}(e^{v_\lambda(x)}-e^{v_\lambda(y)})^2\\
&=\int_V|\nabla e^{v_\lambda}|^2 d\mu.
\end{split}\]
Inserting this estimate into \eqref{eq:eq1}, we obtain
\begin{equation}\label{eq:less}
\int_V|\nabla e^{v_\lambda}|^2 d\mu+\int_V\rho e^{2(v_\lambda+w_1)}d\mu-\int_VK_\lambda e^{2(2v_\lambda+w_1)}d\mu\leq0.
\end{equation}
On the other hand, let $h=e^{v_\lambda}$ in \eqref{eq:convex}, yields
\[\int_V |\nabla e^{v_\lambda}|^2d\mu-2\int_VK_\lambda e^{2(2v_\lambda+w_1)}d\mu\geq0.\]
Together with \eqref{eq:less}, we have
\begin{equation}\label{eq:est}
\int_{V}|\nabla e^{v_\lambda}|^2d\mu\leq -2\int_V\rho e^{2(v_\lambda+w_1)}d\mu=-2\rho(x_0)e^{2u_\lambda(x_0)}.
\end{equation}
One may derive that $u_\lambda(x_0)$ is the uniform bound, i.e. there exists a constant $C>0$ such that
\begin{equation*}\label{eq:upper}
e^{2u_\lambda(x_0)}\leq C.
\end{equation*}
Indeed, if $u_\lambda(x_0)\leq0,$ the above inequality is obvious; if $u_\lambda(x_0)>0,$  as in the proof of Lemma \ref{lemma:ps}, see \eqref{eq:con}, one can get the boundedness of $u_\lambda(x_0)$ as well. Together with \eqref{eq:est}, yields
\begin{equation}\label{eq:gra}\int_{V}|\nabla e^{v_\lambda}|^2d\mu\leq C'.\end{equation}
Next, we \textit{claim} that $e^{v_\lambda}$ and thus $e^{u_\lambda}$ is uniformly bounded in $\ell^2_\mu$.
Suppose not, we may assume that $\|e^{v_\lambda}\|_{\ell^2_\mu}\rightarrow \infty$ as $\lambda\rightarrow \lambda^*.$  Let
\[w_\lambda=\frac{e^{v_\lambda}}{\|e^{v_\lambda}\|_{\ell^2_\mu}},\]
then $\|w_\lambda\|_{\ell^2_\mu}=1$, and $\|\nabla w_\lambda\|_{\ell^2_\mu}\rightarrow 0$ from \eqref{eq:gra}. It follows that $w_\lambda$ converges to a constant in $W^{1,2}(V)$. From \eqref{eq:upper}, we have $w_\lambda(x_0)\rightarrow 0$, and hence $w_\lambda\equiv0$ on $V$, which contradicts $\|w_\lambda\|_{\ell^2_\mu}=1$ and then confirms our claim. Since $u_\lambda$ is bounded below by \eqref{lower} in Step 1, one has the $\ell^2_\mu$-boundedness of $u_\lambda$, and thus the $W^{1,2}(V)$-boundedness of $u_\lambda$. $\hfill\Box$

\subsection{No solution when $\lambda>\lambda^*$.}
Proof of this case is a consequence of  the upper and lower solutions principle, as follows.

\textit{Proof of (4) in Theorem \ref{thm:main}.}
Let $u_{\lambda_1}$ be the solution of \eqref{eq:0} at some $\lambda_1>\lambda^*$. For any $0<\lambda<\lambda_1$, $u_{\lambda_1}$ is an upper solution of \eqref{eq:0} at $\lambda$. Indeed,
\[\Delta u_{\lambda_1}+\kappa-(K+\lambda)e^{2u_{\lambda_1}}=(\lambda_1-\lambda)e^{2u_{\lambda_1}}>0.\]
From \eqref{eq:lo-solu}, it is easy to get a lower bound solution $\varphi$ of \eqref{eq:0} at $\lambda$ such that $\varphi\leq u_{\lambda_1}$. By the upper and lower solutions principle, there exists a solution of \eqref{eq:0} at $\lambda$, which contradicts the definition of $\lambda^*.$
$\hfill\Box$

\hspace{3cm}\\

\end{document}